\documentclass[11pt]{amsart}
\usepackage{amssymb,amsmath,mathrsfs,cite}

\newcommand{\R}{{\mathbb R}}
\newcommand{\Hh}{{\mathbb H}}

\newcommand{\be}{\begin{equation}}
\newcommand{\ee}{\end{equation}}

\newcommand{\eps}{\varepsilon}

\renewcommand{\k}{\kappa}

\numberwithin{equation}{section}
\newtheorem{theorem}{Theorem}[section]
\newtheorem{proposition}[theorem]{Proposition}
\newtheorem{corollary}[theorem]{Corollary}
\newtheorem{lemma}[theorem]{Lemma}

\theoremstyle{definition}
\newtheorem{remark}[theorem]{Remark}

\newcommand{\brm}{\begin{remark}\rm}
\newcommand{\erm}{\end{remark}}
\newcommand{\bte}{\begin{theorem}}
\newcommand{\ete}{\end{theorem}}
\newcommand{\bpr}{\begin{proposition}}
\newcommand{\epr}{\end{proposition}}
\newcommand{\ble}{\begin{lemma}}
\newcommand{\ele}{\end{lemma}}
\newcommand{\beq}{\begin{equation}}
\newcommand{\eeq}{\end{equation}}
\newcommand{\bdm}{\begin{displaymath}}
\newcommand{\edm}{\end{displaymath}}
\numberwithin{equation}{section}

\title[On the spatial segregation limit for a competition-diffusion 
system]{On the long term spatial segregation \\ for a competition-diffusion system}

\author{Marco Squassina}

\address{Dipartimento di Informatica
\newline\indent
Universit\`a degli Studi di Verona
\newline\indent
C\'a Vignal 2, Strada Le Grazie 15
\newline\indent
I-37134 Verona, Italy}
\email{marco.squassina@univr.it}

\thanks{The author was supported by the MIUR national research
project ``Variational and Topological Methods in the Study of
Nonlinear Phenomena''}

\subjclass[2000]{35B40, 35K57, 35B35, 92D25}
\keywords{Competition-diffusion systems, Lotka-Volterra model, spatial segregation, 
population dynamics, asymptotic behaviour, stationary solution, dissipative systems}

\begin{document}
\begin{abstract}
We investigate the long term behavior for a class of competition-diffusion systems
of Lotka-Volterra type for two competing species in the case of low regularity assumptions on the data.
Due to the coupling that we consider the system cannot be reduced to a single equation yielding
uniform estimates with respect to the inter-specific competition rate parameter.
Moreover, in the particular but meaningful case of  initial data with disjoint support and Dirichlet boundary
data which are time-independent, we prove that as the competition rate goes to infinity the solution
converges, along with suitable sequences, to a spatially segregated state satisfying 
some variational inequalities. 
\end{abstract}
\maketitle

%----------------------------------------------------------------------

\section{Introduction}

Let $\Omega$ be a bounded, open, connected subset of $\R^N$ with smooth boundary and
let $\kappa$ be a positive parameter. The aim of this paper is to investigate the 
asymptotic behavior  of a competition-diffusion system of Lotka-Volterra type
for two competing species of population of densities $u$ and $v$, with Dirichlet boundary conditions, 
\begin{equation}
\label{system}
\tag{$P_{\kappa}$}
\begin{cases}
u_t-\Delta u=f(u)-\kappa uv^2, & \text{in $\Omega\times(0,\infty)$}, \\
v_t-\Delta v=g(v)-\kappa vu^2, & \text{in $\Omega\times(0,\infty)$}, \\
u(x,t)=\psi(x,t), & \text{on $\partial\Omega\times[0,\infty)$}, \\
v(x,t)=\zeta(x,t), & \text{on $\partial\Omega\times[0,\infty)$}, \\
u(x,0)=u_0(x), &\text{in $\Omega$}, \\
v(x,0)=v_0(x), &\text{in $\Omega$}.
\end{cases}
\end{equation}
A relevant problem in  population ecology is the understanding of the interactions
between different species, in particular in the case when the interactions are large and of competitive type.
As the inter-specific parameter $\kappa$ ruling the mutual interaction of the species gets large, competitive
reaction-diffusion systems are expected to approach
a limiting configuration where the populations survive but 
exhibit disjoint habitats (cf.\ \cite{shige,nambamim,dancerzang,dancermimura,dancermimura1,mikawa}). 
For population dynamics models which require Dirichlet boundary conditions we refer 
to \cite{nambamim,dancermimura}, while for the more ecologically natural Neumann boundary conditions we refer 
to  \cite{dancerzang,lida} and references therein.  As pointed out in
 \cite{dancermimura}, the Dirichlet case presents further difficulties compared with the 
 Neumann case, as the boundary terms which pop up after integration by parts cannot be estimated  independently of $\kappa$.
The classical stationary Lotka-Volterra model for two populations 
\begin{equation}
\label{orig2po}
\begin{cases}
-\Delta u=f(u)-\kappa uv, & \text{in $\Omega$,} \\
-\Delta v=g(v)-\kappa vu, & \text{in $\Omega$,} \\
u=\psi, & \text{on $\partial\Omega$}, \\
v=\zeta, & \text{on $\partial\Omega$}
\end{cases}
\end{equation}
has been intensively studied with respect to the spatial segregation limit as $\kappa\to\infty$.
If, for instance, $\psi$ and $\zeta$ 
belong to $W^{1,\infty}(\partial\Omega)$, then there exists a sequence of solution $(u_\kappa,v_\kappa)$ to \eqref{orig2po},
bounded in $W^{1,\infty}(\overline{\Omega})$, and a limiting positive state $(u,v)$ 
with $uv=0$, satisfying suitable variational inequalities and such that, 
up to a subsequence, $u_\kappa\to u$ and $v_\kappa\to v$ in $H^1(\Omega)$ with
a precise rate of convergence (see \cite{ctv-adv}). 
Concerning the parabolic system associated with \eqref{orig2po}, in \cite{dancermimura} 
Crooks, Dancer, Hilhorst, Mimura and Ninomiya
proved (also in the case of possibly different diffusion coefficients) that, for any $T>0$, there exists subsequences 
$u_{\kappa_m}$ and $v_{\kappa_m}$ of the solutions converging 
in $L^2(\Omega\times(0,T))$ to a bounded state with disjoint support and solving a limiting free boundary
problem.  Beside this convergence results on finite time intervals, in \cite{dancermimura1}, in the case of equal diffusion coefficients
and stationary boundary conditions,
Crooks, Dancer and Hilhorst recently studied the long term segregation for large interactions,
by reducing the system to a single equation whose solutions admit uniform estimates in $\kappa$.
Typically, stabilization is based upon a variational structure yielding  an energy functional,
bounded  and decreasing along the trajectories (see e.g.\ \cite{haraux,temam}). 
Unfortunately, as far as we know, due to the coupling term $-\kappa uv$,
the parabolic system associated with \eqref{orig2po} does not admit a natural Lyapunov functional 
and a direct analysis is therefore not possible. 
Now, system \eqref{system} can be regarded as a variant of the standard Lotka-Volterra model,
with different inter-specific competition coupling terms. In addition, if one considers
homogeneous boundary data, then \eqref{system} admits a natural non-increasing 
 energy functional $\Lambda_\kappa:[0,\infty)\to\R$
\begin{equation*}
\Lambda_\kappa(t)=\frac 12\|\nabla u(t)\|_2^2+\frac 12\|\nabla v(t)\|_2^2-
\int_\Omega\int_0^{u(t)} f(\sigma)d\sigma
-\int_\Omega\int_0^{v(t)} g(\sigma)d\sigma+
\frac \kappa2\int_\Omega u^2(t)v^2(t).
\end{equation*}
As we will see, a non-increasing energy functional can be constructed also for general
boundary conditions (see the proof of Theorem \ref{stabbz}).
We shall tackle the problem with techniques from the theory of dissipative dynamical systems to
show the convergence towards the solutions to the stationary system, formally written as 
\begin{equation}
\label{limsystem}
\tag{$S_\kappa$}
\begin{cases}
-\Delta u=f(u)-\kappa uv^2, & \text{in $\Omega$,} \\
-\Delta v=g(v)-\kappa vu^2, & \text{in $\Omega$,} \\
u=\psi_\infty, & \text{on $\partial\Omega$}, \\
v=\zeta_\infty, & \text{on $\partial\Omega$}.
\end{cases}
\end{equation}
A question which naturally arises is whether the  solutions  stabilize
towards a  segregated state along some $t_j\to\infty$ and $\kappa_j\to\infty$, for instance in 
the natural case when the initial data have disjoint support and the boundary data are stationary in time 
(see problem \eqref{system-BDtime ind} in the next section). 
Some numerical computations in a square domain in $\R^2$ have been performed in \cite[see Sections 1 and 4]{dancermimura} 
for the Lotka-Volterra model under these assumptions on the initial and boundary conditions (see also \cite{numerikmim},
where an algorithm for parallel computing was implemented in order to efficiently track the interfaces).
In \cite{squazucch} we arranged a complete set of numerical experiments both for \eqref{system-BDtime ind}
(i.e.\ system \eqref{system} with time-independent boundary data) and the corresponding model with the standard Lotka-Volterra coupling.
Although on one hand working with \eqref{system} gives some advantages in
the study of the long term dynamics for $\kappa$ fixed as it directly admits a Lyapunov functional,
on the other hand the  asymptotic analysis  
for the solutions of  \eqref{limsystem} is far more complicated than the study of \eqref{orig2po} 
(subtracting the equations of \eqref{orig2po} one reduces to the single equation $\Delta u=\k u(u-\Phi)$
where $\Phi$ is an harmonic function, while this
is not the case working with \eqref{limsystem}). For instance, the global boundedness in $\k$ of the solutions in $H^1$ 
will be derived from the corresponding boundedness for the solution flow of the parabolic system
uniformly with respect to $\k$. To show the boundedness directly on the elliptic systems seems
out of reach. In addition, the blow up analysis based on Lipschitz rescalings performed in
\cite{ctv-adv} does not seem to work.
\vskip2pt

Concerning some physical motivations to consider coupling terms between the equations which are different from 
the standard one $uv$, we refer the reader, e.g., to Section 3.3 of classical Murray's book \cite{murray}
(looking at formula (3.14) at p.87, our system corresponds to the choice $F(N,P)=1-N-\kappa P^2$ and
$G(N,P)=1-P-\kappa N^2$ with respect to the book's notations).
It is also useful to think about systems of two Schr\"odinger \cite{ablowitz} or Gross-Pitaevskii \cite{dalfovo} equations modelling
particle interaction (and populations can also be thought as discrete collections of interacting particles), intensively investigated 
in recent time (nonlinear optics, Bose-Einstein binary condensates, etc.), 
which present all the coupling of \eqref{system}, yielding a variational 
structure. We refer the reader to \cite{pellmontsquas} 
for the case $\kappa<0$, with physical motivations e.g.\ from \cite{sc-m-neg}, and to \cite{pomponio} for the case where $\kappa>0$, 
with physical motivations e.g.\ from \cite{sc-m-pos}. Both  \cite{pellmontsquas,pomponio} deal with the semi-classical
regime analysis.

\vskip6pt
\noindent
\subsection{The main result}
The main result of the paper concerns with the long-term behaviour
in large-competition regime for the system with time-independent boundary data, that is
\begin{equation}
\label{system-BDtime ind}
\begin{cases}
u_t-\Delta u=f(u)-\kappa uv^2, & \text{in $\Omega\times(0,\infty)$}, \\
v_t-\Delta v=g(v)-\kappa vu^2, & \text{in $\Omega\times(0,\infty)$}, \\
u(x,t)=\psi(x), & \text{on $\partial\Omega\times[0,\infty)$}, \\
v(x,t)=\zeta(x), & \text{on $\partial\Omega\times[0,\infty)$}, \\
u(x,0)=u_0(x), &\text{in $\Omega$}, \\
v(x,0)=v_0(x), &\text{in $\Omega$}.
\end{cases}
\end{equation}
Concerning the functions $f,g:\R\to\R$, let:
\begin{align}
\label{eq:dinint}
& f,g\in C^1([0,\infty)),\quad f(s)=g(s)=0,\quad \text{ for all }s\leq0, \\
\nonumber &f(s)<0,\ g(s)<0,\quad \text{ for all }s>1,
\end{align}
and we set
$$
F(t)=\int_0^t f(\sigma)d\sigma,\qquad
G(t)=\int_0^t g(\sigma)d\sigma.
$$
The initial and boundary data are required to satisfy:
\begin{align}
\label{eq:u0v0-0}
&u_0,v_0\in H^1(\Omega),\quad 0\leq u_0(x)\leq 1,\
0\leq v_0(x)\leq 1, \quad\text{a.e.  in } \Omega, \\
 \noalign{\vskip2pt}
\label{eq:sta-0} & \psi,\zeta\in H^{1/2}(\partial\Omega),\quad
\psi=u_0|_{\partial\Omega},\quad
\zeta=v_0|_{\partial\Omega},\\
 \noalign{\vskip2pt}
\label{eq:psizeta-0}& 0\leq \psi(x)\leq1,\,\,\ 0\leq\zeta(x)\leq1,
\quad\text{on }\partial \Omega.
\end{align}
Under these assumptions, as well as those of Section 2, for all
$\kappa>0$, system \eqref{system} admits a unique global solution $u_\kappa,v_\kappa\in
C^0([0,\infty),H^1(\Omega))\cap C^1((0,\infty),L^2(\Omega))$. For the local existence, we refer 
the reader to a paper by Hoshino-Yamada  \cite{hoshi-yama} (see e.g.\ Theorems 1 and 2, 
having in mind to choose $\theta=\alpha=\gamma=\frac{1}{2}$ in Theorem 1(i) and $\gamma=0$ 
in Theorem 2(ii), with respect to the notations therein).
For smoothing effects we also wish to refer to the classical book of Henry \cite{henry}.
The global existence result can be deduced by the comparison principle for parabolic 
equations (see, for example, the book of Smoller \cite{smoller}).
For $u_t-\Delta u=f(u)-\kappa uv^2$, $v_t-\Delta v=g(v)-\kappa vu^2$
with positive initial data, one can show $0\leq u(t)\leq U(t)$ and $0\leq v(t)\leq V(t)$,
where $U, V$ are solutions of $U_t-\Delta U=f(U), V_t-\Delta V=g(V)$
with the same initial and boundary conditions. Since $U$ and $V$ exist 
globally in time due to assumptions \eqref{eq:dinint}, \eqref{eq:u0v0-0} and \eqref{eq:psizeta}
(a priori uniform-in-time $L^\infty$-estimates for the solutions hold, see Lemma \ref{l:mp}), one also recovers 
the global existence result (for the sake of completeness, we also mention 
Theorem 3 in Hoshino-Yamada  \cite{hoshi-yama} for small initial data
and part (iv) of Proposition 7.3.2 in \cite{lunardi} for smooth initial data). In the following we set
$\Hh=H^1(\Omega)\times H^1(\Omega)$, endowed with the  standard Dirichlet norm, and
$$
\Hh_0=\{(u,v)\in \Hh: \text{$uv=0$ a.e.\ in $\Omega$}\}.
$$ 
\vskip4pt
The following is the main result of the paper, regarding system \eqref{system-BDtime ind}.
\begin{theorem}
\label{main1}
Assume \eqref{eq:dinint}-\eqref{eq:psizeta-0} and $(u_0,v_0)\in\Hh_0$. Let $(u_\kappa,v_\kappa)$
be the solution to system \eqref{system-BDtime ind}. Then there exist two diverging
sequences $(\kappa_m),(t_m)\subset\R^+$ and $(u_\infty,v_\infty)\in \Hh_0$ such that
$$
(u_{\kappa_m}(t_m),v_{\kappa_m}(t_m)) \to (u_\infty,v_\infty)\quad \text{in the $L^p\times L^p$ norm,
for any $p\in[2,\infty)$,}
$$
as $m\to\infty$, where  
$$
u_\infty,v_\infty\geq 0,\quad 
-\Delta u_\infty\leq f(u_\infty),\quad 
-\Delta v_\infty\leq g(v_\infty),\quad
u_\infty|_{\partial\Omega}=\psi,\quad 
v_\infty|_{\partial\Omega}=\zeta.
$$ 
Moreover, in the one-dimensional case, we have
$$
\|(u_{\kappa_m}(t_m),v_{\kappa_m}(t_m)) -(u_\infty,v_\infty)\|_{L^\infty\times L^\infty}\to 0,\quad
\text{as $m\to\infty$}.
$$
\end{theorem}

\vskip8pt

Hence, starting with segregated data, 
the system evolves towards a limiting segregated state satisfying suitable variational inequalities.
As we have previously pointed out, in Sections 1,4 of \cite{dancermimura}, the reader can find
very nice pictures reproducing (for the classical model) these kind of separation phenomena.
Notice that, due to the nonstandard coupling in system \eqref{limsystem} the $H^1$ convergence
seems pretty hard to obtain either working directly on the system (which would require precise quantitative
estimate of the rate of convergence of the solutions to $u_\infty$ and $v_\infty$) or using indirect arguments such
combining blow up analysis with Liouville theorems (which, however, would naturally require 
stronger regularity assumptions on the boundary conditions). In Section \ref{sectev}, we will obtain,
for $\k$ fixed, the asymptotic behaviour of the system in the case of almost stationary boundary data.
The author is not aware of any other result of this type in the literature (see also \cite{chilljend}).

\vskip6pt

\section{Long term behaviour for $\kappa$ fixed}
\label{sectev}

The goal of this section is the study of the long 
term behaviour  of the parabolic system \eqref{system}, for any $\kappa>0$ fixed.
We cover the general case of boundary data  
depending on time. Finally, in the particular case of segregated initial data  
and time independent boundary conditions, we will prove a stronger global
boundedness result.

\subsection{Assumptions  and main result}
Concerning $f$ and $g$  we will assume condition \eqref{eq:dinint}.
Moreover, the initial and boundary data are required to satisfy \eqref{eq:u0v0-0} and
\begin{align}
\label{eq:sta} &\psi,\zeta\in C^0([0,\infty), H^{1/2}(\partial\Omega)),\quad
\psi(0)=u_0|_{\partial\Omega},\quad
\zeta(0)=v_0|_{\partial\Omega},\\
 \noalign{\vskip2pt}
\label{eq:psizeta}& 0\leq \psi(x,t)\leq1,\,\ 0\leq\zeta(x,t)\leq1,
\quad\text{ on }\partial \Omega\times[0,\infty).
\end{align}
We will assume that:
\begin{align}
  \label{eq:2}
  & \psi (\cdot,t)\to\psi_\infty\quad\text{and}\quad
  \zeta (\cdot,t)\to\zeta_\infty\quad \text{in
  }H^{1/2}(\partial\Omega)\quad\text{as }t\to\infty,\\
  \noalign{\vskip4pt}
\label{eq:1}
  & \psi_t,\zeta_t\in  L^1(0,\infty;H^{1/2}(\partial\Omega))\cap  L^2(0,\infty;H^{-1/2}(\partial\Omega)),\,\,\,\,
  \text{$\psi_t (\cdot,t),\zeta_t (\cdot,t)\to 0$ as $t\to\infty$}, \\
   \noalign{\vskip4pt}
  \label{compatib}
&  \psi_t(\cdot,0)=\zeta_t(\cdot,0)=0, \\
\noalign{\vskip4pt}
\label{eq:6}
&\psi_{tt}, \zeta_{tt} \in L^1(0,\infty;H^{-1/2}(\partial\Omega)).
\end{align}
\vskip4pt

Under the previous assumptions we have  the following result.

\begin{theorem}
\label{cve}
Let $(u_0,v_0)\in\Hh$ and $\kappa>0$. Then for every diverging sequence $(t_h)\subset\R^+$ there exist
a subsequence $(t_j)\subset\R^+$ and  a solution $(\hat u_\kappa,\hat v_\kappa)\in\Hh$ to system \eqref{limsystem}
such that
$$
\|(u_\kappa(t_j),v_\kappa(t_j))-(\hat u_\kappa,\hat v_\kappa)\|_{\Hh}\to 0,
\quad\text{as $j\to\infty$}.
$$
Moreover, the convergence holds in the $L^p\times L^p$ norm for any $p\in[2,\infty)$.
\end{theorem}

\vskip4pt
Strenghtening the assumptions we obtain the global boundedness uniformly in $\k$.

\begin{theorem}
\label{cve-cor}
Assume that $(u_0,v_0)\in\Hh_0$ and the boundary conditions are time-independent.
Then, in addition to the conclusion of Theorem \ref{cve}, we have
$$
\sup_{t\geq 0}\sup_{\kappa>0}\|(u_\kappa(t),v_\kappa(t))\|_\Hh<\infty,
$$
namely $(u_\kappa,v_\kappa)$  is  bounded in $\Hh$ (and in any $L^p\times L^p$ space),
uniformly with respect to $\kappa$.
\end{theorem}
This second achievement will be of course an important step in order to
prove the main result of the paper. 

\medskip

\subsection{Some Preliminary results}

From a  direct computation, we have positivity and a priori
bounds for the solutions to \eqref{system}, uniformly with respect to $\kappa$.

\begin{lemma}
\label{l:mp}
$\Sigma=[0,1]\times [0,1]$ is a globally positively invariant region for system \eqref{system},
uniformly with respect to $\k$, namely
$$
0\leq u_\kappa(x,t)\leq 1,\qquad
0\leq v_\kappa(x,t)\leq 1,
\quad\quad\text{a.e. }x\in \Omega,\,\, t\geq 0.
$$
\end{lemma}
\begin{proof}
Testing the first equation of \eqref{system} with
$-u_\kappa^-$ and using \eqref{eq:dinint}, \eqref{eq:u0v0-0}, \eqref{eq:sta} and \eqref{eq:psizeta}, we easily
obtain that $u_\kappa\geq 0$, while testing the same equation
with $(u_{\kappa} -1)^+$  we deduce similarly
that $u_\kappa\leq 1$. An analogous manipulation of the second
equation in \eqref{system} yields the corresponding bounds for the
component $v_\kappa$.
\end{proof}

Let $A=-\Delta$ be the Laplace operator on $L^2(\Omega)$
with domain ${\mathscr D}(A)=H^1_0(\Omega)\cap H^2(\Omega)$ and consider
the hierarchy of Hilbert spaces $H^\alpha={\mathscr D}(A^{\alpha/2})$, $\alpha\in\R$,
with $\|u\|_{H^\alpha}=\|A^{\alpha/2}u\|_{2}$. We recall an exponential decay 
property of the heat kernel operator $e^{t\Delta}$.

\begin{lemma}
\label{expboundlap}
Let $\alpha> 0$. Then there exist $\omega>0$ 
and $C_{\alpha}>0$ such that
\begin{equation}
\label{eq:expoest}
\|e^{t\Delta}\|_{{\mathcal L}(L^2,\,H^{2\alpha})}
\leq C_{\alpha} e^{-\omega t}t^{-\alpha},
\qquad t>0.
\end{equation}
In particular
$$
\int_0^{\infty}\|e^{\sigma\Delta}\|_{{\mathcal L}(L^2,\,H^{2\alpha})}d\sigma<\infty,
$$
provided that $\alpha\in(0,1)$.
\end{lemma}
\begin{proof}
As the real part of the spectrum of $A$ is bounded away from zero 
by a positive constant $\omega$, by \cite[Theorem 1.4.3, p.26]{henry}, 
for $\alpha>0$ there exists $C_\alpha>0$ such that
$\|A^\alpha e^{-tA}\|_{{\mathcal L}(L^2,\,L^2)}\leq C_\alpha e^{-\omega t}t^{-\alpha}$, 
for all $t>0$. Hence,
$\|e^{t\Delta}\|_{{\mathcal L}(L^2,\,H^{2\alpha})}
=\|(-\Delta)^{\alpha} e^{t\Delta}\|_{{\mathcal L}(L^2,\,L^2)}
\leq C_\alpha e^{-\omega t}t^{-\alpha}$, for all $t>0$.
The second assertion follows by \eqref{eq:expoest}.
\end{proof}

Next we provide a compactness result for the trajectories of \eqref{system}.

\begin{lemma}
\label{compemb}
For any $(u_0,v_0)\in \Hh$, $\kappa>0$ and $\tau>0$ the set
$\{(u_\kappa(t),v_\kappa(t)):t\geq\tau\}$ is relatively compact in $\Hh$.
\end{lemma}
\begin{proof}
Let $U$ and $V$ denote the solutions to the linear problems 
\begin{equation}
\label{linearU}
\begin{cases}
U_t-\Delta U=0, &\text{in }\Omega\times(0,\infty),\\
\noalign{\vskip2pt}
U(x,t)=\psi(x,t) , &\text{on }\partial\Omega\times(0,\infty),\\
\noalign{\vskip2pt}
  U(x,0)=U_0(x), & \text{in $\Omega$},
\end{cases}
\end{equation}
and
\begin{equation}
\label{linearV}
\begin{cases}
V_t-\Delta V=0,&\text{in }\Omega\times(0,\infty),\\
\noalign{\vskip2pt}
V(x,t)=\zeta(x,t) , &\text{on }\partial\Omega\times(0,\infty),\\
\noalign{\vskip2pt}
  V(x,0)=V_0(x), & \text{in $\Omega$},
\end{cases}
\end{equation}
where $U_0,V_0\in H^1(\Omega)$ satisfy
$$
\begin{cases}
-\Delta U_0=0,&\text{in }\Omega,\\
\noalign{\vskip2pt}
\,\, U_0(x)=\psi(x,0),&\text{on }\partial\Omega,
\end{cases}
\quad \qquad
\begin{cases}
-\Delta V_0=0,&\text{in }\Omega,\\
\noalign{\vskip2pt}
\,\, V_0(x)=\zeta(x,0),&\text{on }\partial\Omega.
\end{cases}
$$
By assumption \eqref{eq:psizeta} and the maximum principle
for harmonic functions, $0\leq U_0(x)\leq 1$ and $0\leq V_0(x)\leq 1$ for a.e.\ $x\in\Omega$.
Hence, arguing as in the proof of Lemma \ref{l:mp}, we have
$0\leq U(x,t)\leq 1$ and $0\leq V(x,t)\leq 1$ for a.e.\ $x\in \Omega$ and $t\geq 0$.
Now, the functions 
\begin{equation}
\label{link}
\tilde u_\kappa(x,t)= u_\kappa(x,t)-U(x,t),\qquad
\tilde v_\kappa(x,t)= v_\kappa(x,t)-V(x,t)
\end{equation}
solve the system with homogeneous  boundary conditions
\begin{equation}
\label{tildesystem}
\tag{$\widetilde P_{\kappa}$}
\begin{cases}
  (\tilde u_\kappa)_t-\Delta \tilde u_\kappa =f(\tilde u_\kappa+U)-\kappa (\tilde u_\kappa+U)(\tilde v_\kappa+V)^2, & \text{in $\Omega\times(0,\infty)$}, \\
  \noalign{\vskip2pt}
  (\tilde v_\kappa)_t-\Delta \tilde v_\kappa =g(\tilde v_\kappa+V)-\kappa (\tilde v_\kappa+V)(\tilde u_\kappa+U)^2, & \text{in $\Omega\times(0,\infty)$}, \\
   \noalign{\vskip2pt}
  \,\tilde u_\kappa(x,t)=\tilde v_\kappa(x,t)=0, & \text{on $\partial\Omega\times[0,\infty)$}, \\
   \noalign{\vskip2pt}
  \tilde u_\kappa(x,0)=u_0(x)-U_0(x), &\text{in $\Omega$}, \\
   \noalign{\vskip2pt}
  \tilde v_\kappa(x,0)=v_0(x)-V_0(x), &\text{in $\Omega$}.
\end{cases}
\end{equation}
Denote now by $\Psi=\Psi(x;t)\in  C^0([0,\infty),H^{1}(\Omega))$  the family of harmonic extensions to
$\Omega$ of $\psi$
\begin{equation}
\label{harm-ext-psi}
\begin{cases}
-\Delta \Psi(x;t)=0,&\text{in }\Omega,\\
\noalign{\vskip2pt}
\,\, \Psi(x;t)=\psi(x,t), &\text{on }\partial\Omega,
\end{cases}
\end{equation}
and set $\bar U(x,t)=U(x,t)-\Psi(x;t)$. Then  $\bar U$ solves the 
nonautonomous problem with homogeneous boundary and initial conditions
\begin{equation}
\label{0bb}
\begin{cases}
  \bar U_t-\Delta \bar U=-\Psi_t, & \text{in
    $\Omega\times(0,\infty)$},\\
  \bar U(x,t)=0, & \text{on $\partial\Omega\times(0,\infty)$},\\
  \bar U(x,0)=0, &\text{in $\Omega$}.
\end{cases}
\end{equation}
Notice that $ \bar U(x,0)=0$ since $U_0(x)$ and $\Psi(x;0)$ are both harmonic functions
with the same boundary conditions. 
From \eqref{compatib}-\eqref{eq:6} and classical regularity theory for harmonic functions,   
\begin{equation*}
\|\Psi_t\|_{L^\infty(0,\infty;L^2(\Omega))}
\leq c\|\psi_t\|_{L^\infty(0,\infty;H^{-1/2}(\partial\Omega))}
\leq \|\psi_{tt}\|_{L^1(0,\infty;H^{-1/2}(\partial\Omega))}.
\end{equation*}
By Duhamel's formula $\bar U$ is given by
\begin{equation*}
\bar U(t)= -\int_0^t e^{(t-\sigma)\Delta}\Psi_t(\sigma)\,d\sigma.
\end{equation*}
If $\alpha\in(1/2,1)$, in light of \eqref{eq:expoest} of Lemma \ref{expboundlap}, 
since $\Psi_t$ is in $L^\infty(0,\infty;L^2(\Omega))$,
\begin{equation}
\label{eq:est1comp}
\sup_{t\geq 0}\|\bar U(t)\|_{H^{2\alpha}}<\infty.
\end{equation}
Of course the same control holds for $\bar V(t)$.
Let now $\Psi_\infty$ denote the harmonic extension of $\psi_\infty$,
the limit of $\psi(t)$ in $H^{1/2}(\partial\Omega)$ as $t\to\infty$ according to \eqref{eq:2} . 
By standard regularity estimates,
$\|\Psi(t)-\Psi_\infty\|_{H^1(\Omega)}\leq c\|\psi(t)-\psi_\infty\|_{H^{1/2}(\partial\Omega)}$, 
so that $\Psi(t)\to \Psi_\infty$ in $H^1(\Omega)$ as $t\to\infty$. Of course the same control 
holds for the boundary extensions of $\zeta$. Also, by Duhamel's formula we have
\begin{align*}
\tilde u_\kappa(t) &= e^{t\Delta}(u_0-U_0) + \int_0^t e^{(t-\sigma)\Delta}\Phi^1_\kappa(\sigma) d\sigma,  \\
\tilde v_\kappa(t)  &= e^{t\Delta}(v_0-V_0) + \int_0^t e^{(t-\sigma)\Delta}\Phi^2_\kappa(\sigma) d\sigma,
\end{align*}
where 
\begin{equation*}
\Phi^1_\kappa(\sigma) 
=f(u_\kappa(\sigma))-\kappa u_\kappa(\sigma)v_\kappa^2(\sigma),  \qquad
\Phi^2_\kappa(\sigma) 
=g(v_\kappa(\sigma))-\kappa v_\kappa(\sigma)u_\kappa^2(\sigma).
\end{equation*}
By means of Lemma \ref{l:mp}, we have $\Phi^1_\kappa,\,\,\Phi^2_\kappa\in L^\infty(0,\infty;L^\infty(\Omega))$.
If $\alpha\in(1/2,1)$, then again by \eqref{eq:expoest} one obtains for any $\tau>0$
\begin{equation}
\label{eq:est1comp2}
\sup_{t\geq \tau}\|\tilde u_\kappa(t)\|_{H^{2\alpha}}<\infty.
\end{equation}
As $H^{2\alpha}$ is compactly embedded in $H^1(\Omega)$ and
$u_\kappa(t)=\tilde u_\kappa(t)+\bar U(t)+\Psi(t)$
the assertion follows by \eqref{eq:est1comp}-\eqref{eq:est1comp2} for the
component $u_\kappa$. The same arguments works for $\tilde v_\kappa$.
\end{proof}

\begin{remark}
By strengthening the regularity assumptions on the boundary data, say $W^{1,\infty}(\partial\Omega)$
in place of $H^{1/2}(\partial\Omega)$ in the assumptions at the beginning of the section, and defining $-\Delta$ over $L^q(\Omega)$ for any $q\geq 2$,
the previous result can of course be improved, yielding compactness of the trajectories in 
$W^{2\alpha,q}(\Omega)$ for any $q\geq 2$, and hence into spaces of H\"older continuous functions.
Unfortunately the estimates  are not independent
of $\k$ and in order to have $H^1$ bounds uniformly in $\k$ we shall need
to exploit energy arguments. 
\end{remark}

\noindent For every $\tau>0$ and every function $h:(0,\infty)\to
H^1(\Omega)$, let us set
$$
h^\tau(t)=h(t+\tau),\qquad t>0.
$$
The following result gives a stabilization property for the solutions of
the linear parabolic equation with nonhomogeneous time-dependent
boundary conditions.

\begin{lemma}
\label{stabilU}
Let $U$ be the solution to the problem \eqref{linearU}.
Then $U(t)\to U_{\infty}$ in $H^1(\Omega)$ as $t\to\infty$, where
$U_{\infty}\in H^1(\Omega)$ is the solution to 
\begin{equation}
\label{eq:uinfty}
\begin{cases}
-\Delta U_{\infty}=0,&\text{in }\Omega,\\
U_{\infty}=\psi_{\infty},&\text{on }\partial\Omega.
\end{cases}
\end{equation}
\end{lemma}
\begin{proof}
With the notations introduced in the proof of Lemma \ref{compemb}, we consider, for $\tau>0$,
the functions $W(t)=\bar U^\tau(t)-\bar U(t)$ and $\varrho(t)=\Psi_t(t)-\Psi_t^\tau(t)$, which satisfy
\begin{equation}
\label{eq:w}
\begin{cases}
W_t-\Delta W=\varrho(t), & \text{in $\Omega\times(0,\infty)$},\\
W(x,t)=0, & \text{on $\partial\Omega\times(0,\infty)$},\\
W(x,0)=U(\tau)-U_0+\Psi(0)-\Psi(\tau), & \text{in $\Omega$}.
\end{cases}
\end{equation}
By multiplying the equation by $-\Delta W$, we get
$$
\frac{d}{dt}\|\nabla W(t)\|_{2}^2+\|\Delta
W(t)\|_{2}^2=-\int_\Omega \varrho(t)\Delta W(t).
$$
By applying H\"older and then Young inequalities on the right-hand side, we have
$$
\frac{d}{dt}\|\nabla W(t)\|_{2}^2+\frac{1}{2}\|\Delta
W(t)\|_{2}^2\leq \frac{1}{2}\|\varrho(t)\|_2^2.
$$
Let $A$ be the positive operator on $L^2(\Omega)$
defined by $A=-\Delta$, with domain ${\mathcal D}(A)=H^2(\Omega)\cap H^1_0(\Omega)$.
Due to the (compact and dense) injection $H^2(\Omega)\cap H^1_0(\Omega)
={\mathcal D}(A)\hookrightarrow {\mathcal D}(A^{1/2})=H^1_0(\Omega)$,
we have $\|\nabla W\|_2\leq \alpha_1^{-1/2}\|\Delta W\|_2$ for some $\alpha_1>0$ 
(see e.g.\  Henry \cite{henry}), so that
$$
\frac{d}{dt}\|\nabla W(t)\|_{2}^2+\frac{\alpha_1}{2}\|\nabla
W(t)\|_{2}^2\leq \frac{1}{2}\|\varrho(t)\|_2^2.
$$
Finally, Gronwall inequality entails
$$
\|W(t)\|_{H^1_0}^2\leq
\|W(0)\|_{H^1_0}^2e^{-\sigma t}+ce^{-\sigma
t}\int_0^te^{\sigma s}\|\varrho(s)\|_{2}^2ds,
$$
for some $\sigma>0$ and $c>0$. In turn, we readily obtain
\begin{equation*}
\lim_{t\to\infty}\|U^\tau(t)-U(t)\|_{H^1} \leq
  \frac{c}{\sqrt{\sigma}}\lim_{t\to\infty}\|\Psi_t^\tau(t)-\Psi_t(t)\|_{2}
 +\lim_{t\to\infty}\|\Psi^\tau(t)-\Psi(t)\|_{H^1}.
\end{equation*}
In view of \eqref{eq:2}  and standard elliptic equations, we deduce 
\begin{equation}
\label{convrgU}
\|U^\tau(t)-U(t)\|_{H^1}\to 0,\qquad\text{as $t\to\infty$}. 
\end{equation}
The same argument 
shows that $\{U(t)\}_{t\geq 0}$ is bounded in $H^1(\Omega)$. Let now $(t_h)\subset\R^+$ be any diverging sequence.
Since $\{U(t)\}_{t\geq 1}$ is relatively compact in
$H^1(\Omega)$, there exists a subsequence, that we still denote by $(t_h)$, such that
$U(t_h)\to U_{\infty}$ in $H^1(\Omega)$. Let $\eta\in C^\infty_c(\Omega)$. By integrating the equation for
$U$ on $(t_h,t_h+1)\times\Omega$, yields
$$
\lim_h\left[\int_{t_h}^{t_h+1}\int_\Omega U_t\eta+
\int_{t_h}^{t_h+1}\int_\Omega \nabla U\cdot\nabla\eta\right]=0.
$$
On one hand, we have
$$
\lim_{h}\left|\int_{t_h}^{t_h+1}\int_\Omega
U_t\eta\right|\leq \lim_{h}\int_\Omega
|U(t_h+1)-U(t_h)||\eta| \leq
c\lim_h\|U^1(t_h)-U(t_h)\|_{2}=0.
$$
Moreover, there exists  $(s_h)\subset\R^+$ with
$s_h=t_h+\xi_h$, $0\leq \xi_h\leq 1$, such that by \eqref{convrgU}
$$
\int_{t_h}^{t_h+1}\!\!\!\int_\Omega\nabla U\cdot\nabla\eta
=\int_\Omega\nabla U(s_h)\cdot\nabla\eta=\int_\Omega\nabla
U(t_h)\cdot\nabla\eta+o(1),\qquad \text{as $h\to\infty$}.
$$
Hence, taking the limit as $h\to\infty$, we get $\int_\Omega\nabla U_{\infty}\cdot\nabla\eta=0$.
Moreover, from the convergence of $U(t_h)$ to $U_{\infty}$ in $H^1(\Omega)$ we deduce
that $U(t_h)|_{\partial\Omega}\to U_{\infty}|_{\partial\Omega}$ in
$H^{1/2}(\partial\Omega)$. From \eqref{eq:2} we deduce that
$U_{\infty} =\psi_{\infty}$ on $\partial\Omega$.  Therefore $U_{\infty}$ solves \eqref{eq:uinfty}.
Since \eqref{eq:uinfty} has a unique solution, we
actually deduce the convergence of the whole flow $U(t)$.
\end{proof}

Next, we obtain a summability result for the solutions to  \eqref{linearU}.

\begin{lemma}
\label{stabilU1}
Let $U$ be the solution to  \eqref{linearU}. Then $U_t\in L^1(0,\infty;H^{1}(\Omega))$.
\end{lemma}
\begin{proof}
As in the proof of Lemma \ref{compemb}, $\bar U$ is the solution to \eqref{0bb}.
Hence, taking into account \eqref{compatib}, it turns out that $\widetilde U(x,t)=
\bar U_t(x,t)$ is a solution to 
\begin{equation}
\label{Vbb}
\begin{cases}
  \widetilde U_t-\Delta \widetilde U=-\Psi_{tt}, & \text{in
    $\Omega\times(0,\infty)$},\\
  \widetilde U(x,t)=0, & \text{on $\partial\Omega\times(0,\infty)$},\\
  \widetilde U(x,0)=0, &\text{in $\Omega$}.
\end{cases}
\end{equation}
By assumption \eqref{eq:6} it follows $\Psi_{tt}\in L^1(0,\infty;L^2(\Omega))$.
In addition, we have $\Psi_{t}\in L^1(0,\infty;H^1(\Omega))$.
By Lemma \ref{expboundlap} we have
$\|e^{t\Delta}\|_{{\mathcal L}(L^2(\Omega),\,H^1_0(\Omega))}\leq Ce^{-\omega t}t^{-1/2}$,
for some $C,\omega>0$. Hence,
\begin{equation*}
\widetilde U(t)= -\int_0^t e^{(t-\sigma)\Delta}\Psi_{tt}(\sigma)\,d\sigma,
\end{equation*}
and we obtain
\begin{align*}
\|\widetilde U\|_{L^1(0,\infty;H^1_0(\Omega))}  &\leq
C\int_0^{\infty}\bigg[\int_0^t e^{-\omega(t-\sigma)}(t-\sigma)^{-1/2}
\|\Psi_{tt}(\sigma)\|_{2}d\sigma\bigg]dt\\
& =C\int_0^{\infty}\|\Psi_{tt}(\sigma)\|_{2}\bigg[\int_\sigma^{\infty}
e^{-\omega(t-\sigma)}(t-\sigma)^{-1/2}\,dt\bigg]d\sigma \\
&=
C\bigg(\int_0^{\infty}
e^{-\omega \sigma}\sigma^{-1/2}\,d\sigma\bigg)\|\Psi_{tt}\|_{L^1(0,\infty;L^2(\Omega))}.
\end{align*}
Hence $\widetilde U\in L^1(0,\infty;H^{1}_0(\Omega))$, 
which yields $\bar U_t\in L^1(0,\infty;H^{1}_0(\Omega))$ and,
in turn, taking into account  \eqref{eq:1}, also $U_t\in L^1(0,\infty;H^{1}(\Omega))$,
concluding the proof.
\end{proof}

\begin{lemma}
\label{integ2summab}
Let $\tilde u_\kappa$ and $\tilde v_\kappa$ be as in system $\widetilde P_{\kappa}$. Then
$$
\int_0^T\|\partial_t \tilde u_\kappa(\sigma)\|_2^2\,d\sigma<\infty,\quad
\int_0^T\|\partial_t \tilde v_\kappa(\sigma)\|_2^2\,d\sigma<\infty,
$$
for any $T>0$.
\end{lemma}
\begin{proof}
Setting $\Upsilon(x,t)=f(\tilde u_\kappa(x,t)+U(x,t))-\kappa (\tilde u_\kappa(x,t)+U(x,t))(\tilde v_\kappa(x,t)+V(x,t))^2$
for any $x\in\Omega$ and $t>0$ and $m(x)=u_0(x)-U_0(x)$, it follows that $\tilde u_\kappa$ is the solution to 
\begin{equation*}
\begin{cases}
  \partial_t\tilde u_\kappa-\Delta \tilde u_\kappa =\Upsilon & \text{in $\Omega\times(0,\infty)$}, \\
  \noalign{\vskip2pt}
  \,\tilde u_\kappa(x,t)=0, & \text{on $\partial\Omega\times[0,\infty)$}, \\
   \noalign{\vskip2pt}
  \tilde u_\kappa(x,0)=m(x), &\text{in $\Omega$}.
\end{cases}
\end{equation*}
Hence, since $m\in H^1_0(\Omega)$ and $\Upsilon\in L^2(0,T,L^2(\Omega))$
for any $T>0$ (as $0\leq  u_\kappa, v_\kappa\leq 1$ and $f$ is continuous), 
the desired summability for $\partial_t\tilde u_\kappa$ follows, e.g., by \cite[Theorem 5, p.360]{evans}.
The proof for $\partial_t\tilde v_\kappa$ is similar.
\end{proof}

Let us recall a useful elementary Gronwall type inequality.
\begin{lemma}
\label{gronwalltype}
Let $g\in L^1([0,\infty),[0,\infty))$. Assume that $\Upsilon:[0,\infty)\to[0,\infty)$ is an 
absolutely continuous function such that
$$
\Upsilon(t)\leq c_1+c_2\int_0^t g(\sigma)\sqrt{\Upsilon(\sigma)}d\sigma,\qquad t\geq 0,
$$
for some $c_1,c_2>0$. Then
$$
\Upsilon(t)\leq 2c_1+c_2^2\|g\|_{L^1(0,\infty)}^2,\qquad t\geq 0.
$$
\end{lemma}
\begin{proof}
Let $t>0$ and consider $\bar t\in [0,t]$ such that $\Upsilon(\bar t)=\max\{\Upsilon(\sigma):\sigma\in[0,t]\}$. Hence
$$
\Upsilon(\bar t)\leq c_1+c_2\int_0^{\bar t} g(\sigma)\sqrt{\Upsilon(\sigma)}d\sigma
\leq c_1+c_2\sqrt{\Upsilon(\bar t)}\|g\|_{L^1(0,t)}
\leq c_1+c_2\sqrt{\Upsilon(\bar t)}\|g\|_{L^1(0,\infty)},
$$
so the assertion immediately follows by Young inequality and $\Upsilon (t)\leq \Upsilon (\bar t)$.
\end{proof}

Next we obtain an $H^1$ stabilization result  for the solutions $(u_\kappa,v_\kappa)$ to \eqref{system}.

\begin{theorem}
\label{stabbz}
Assume that $(u_0,v_0)\in\Hh$ and set
\begin{equation}
 \label{cruquant}
\mu=\|u_0v_0\|_2^2+\|\Psi_t\|_{L^1(0,\infty;L^2(\Omega))}+\|Z_t\|_{L^1(0,\infty;L^2(\Omega))}.
\end{equation}
Then there exists a positive constant $R=R(u_0,v_0,\psi,\zeta)$ independent of $\kappa$ such that
\begin{equation}
 \label{acca1}
\|(u_\kappa(t), v_\kappa(t))\|_\Hh\leq R+\kappa\mu,
\qquad\text{for all $t\geq 0$}.
 \end{equation}
Moreover, for any $\tau_0>0$ and $\kappa>0$,
$$
\lim_{t\to\infty}\sup_{\tau\in[0,\tau_0]}\|u_\kappa(t+\tau)-u_\kappa(t)\|_{H^1}=0,\qquad
\lim_{t\to\infty}\sup_{\tau\in[0,\tau_0]}\|v_\kappa(t+\tau)-v_\kappa(t)\|_{H^1}=0.
$$
\end{theorem}

\begin{proof}
Let $\tau_0>0$ and $\kappa>0$. Let us first prove that
\begin{equation}
\label{firstst}
\lim_{t\to\infty}\sup_{\tau\in[0,\tau_0]}\|u_\kappa(t+\tau)-u_\kappa(t)\|_{2}=0,\qquad
\lim_{t\to\infty}\sup_{\tau\in[0,\tau_0]}\|v_\kappa(t+\tau)-v_\kappa(t)\|_{2}=0.
\end{equation}
According to the proof of Lemma \ref{compemb}, let again $U$ (resp.\ $V$) be the solution of the  
linear problems \eqref{linearU} (resp.\ \eqref{linearV}),
where $U_0$ (resp.\ $V_0$) is the harmonic extensions of $\psi(0)$
(resp.\ $\zeta(0)$). Then $\tilde u_\kappa(x,t)= u_\kappa(x,t)-U(x,t)$ and
$\tilde v_\kappa(x,t)= v_\kappa(x,t)-V(x,t)$ are solutions to system \eqref{tildesystem}
having homogeneous boundary conditions. Let now $\eps\in(0,1)$ 
and, taking into account Lemma \ref{integ2summab}, introduce the auxiliary energy functional $\Lambda_\kappa:[0,\infty)\to\R$ 
defined by setting:
\begin{align*}
  &\Lambda_\kappa(t)=\frac 12\|\nabla \tilde u_\kappa(t)\|_2^2 +\frac
  12\|\nabla \tilde v_\kappa(t)\|_2^2-
  \int_\Omega F(\tilde u_\kappa(t)+U(t))\\
  &\quad-\int_\Omega G(\tilde v_\kappa(t)+V(t))+
  \frac \kappa2\int_\Omega (\tilde u_\kappa(t)+U(t))^2(\tilde v_\kappa(t)+V(t))^2\\
   &\quad+2\int_0^t\bigg[\int_{\Omega} \nabla \tilde u_\kappa(\sigma)\cdot\nabla
   U_t(\sigma)\bigg]d\sigma -
 \int_\Omega \nabla U(t)\cdot\nabla \tilde
    u_\kappa(t)-\int_0^t 
{}_{H^{-\frac 12}\!\!}\bigg\langle
   \frac{\partial \tilde u_\kappa(\sigma)}{\partial\nu},\psi_t(\sigma)
 \bigg\rangle_{\!\!H^{\frac 12}}\!d\sigma\\
&\quad+2\int_0^t\bigg[\int_\Omega \nabla \tilde v_\kappa(\sigma)\cdot\nabla
  V_t(\sigma)\bigg]d\sigma -\int_\Omega \nabla V(t)\cdot\nabla \tilde
  v_\kappa(t)-\int_0^t {}_{H^{-\frac 12}\!\!}\bigg\langle
  \frac{\partial \tilde v_\kappa(\sigma)}{\partial\nu},\zeta_t(\sigma)
\bigg\rangle_{\!\!H^{\frac 12}}\!d\sigma  \\
&\quad
+\eps\int_0^t\|\partial_t \tilde u_\kappa(\sigma)\|_2^2\,d\sigma
+\eps\int_0^t\|\partial_t \tilde v_\kappa(\sigma)\|_2^2\,d\sigma.
\end{align*}
We prove that $\Lambda_\kappa$ is nonincreasing and 
there exist two constants $\alpha_\kappa\in\R$ and $\beta_\kappa\in\R$ (which we
will write down explicitely) such that
$\alpha_\kappa\leq\Lambda_\kappa(t)\leq\beta_\kappa$, for all $t\geq 0$.
By multiplying the first equation of \eqref{tildesystem} by $\partial_t u_\kappa$ and the second
one by $\partial_t v_\kappa$, using the fact that $U$ and $V$ solve problems \eqref{linearU}-\eqref{linearV},
and adding the resulting identities, we reaches
\begin{equation}
\label{ident}
\frac{d}{dt}\Lambda_\kappa(t)=
-(1-\eps)\|\partial_t \tilde u_\kappa(t)\|_2^2
-(1-\eps)\|\partial_t \tilde v_\kappa(t)\|_2^2\leq 0.
\end{equation}
In particular $\{t\mapsto\Lambda_\kappa(t)\}$ is a nonincreasing function. Hence, 
\begin{align*}
\Lambda_\kappa(t)\leq\Lambda_\kappa(0)&=\frac 12\|\nabla (u_0-U_0))\|_2^2
+\frac 12\|\nabla (v_0-V_0)\|_2^2-
\int_\Omega F(u_0)-\int_\Omega G(v_0)\\
&- \int_\Omega \nabla U_0\cdot\nabla(u_0-U_0)
-\int_\Omega \nabla V_0\cdot\nabla(v_0-V_0)
+\frac \kappa2\int_\Omega u_0^2v_0^2,
\end{align*}
for all $t\geq 0$, namely $\Lambda_\kappa$ is bounded from above, uniformly in
time and $\beta_\kappa$ is of the form
\begin{equation}
\label{constantupident}
\beta_\kappa=P+\kappa\|u_0v_0\|_2^2,\qquad P=P(u_0,v_0,\psi,\zeta).
\end{equation}
Now, using the trace inequality, the first equation of $(\widetilde P_\kappa)$,
the $L^\infty$-boundedness of the solutions and the Young inequality, we find $c>0$ and $c_\eps>0$ such that
\begin{align*}
\left|  \int_0^t 
{}_{H^{-\frac 12}\!\!}\bigg\langle
   \frac{\partial \tilde u_\kappa(\sigma)}{\partial\nu},\psi_t(\sigma)
 \bigg\rangle_{\!\!H^{\frac 12}}\!\! d\sigma\right| &\leq 
 \int_0^t \|\nabla\tilde u_\kappa(\sigma)\|_2\|\nabla\Psi_t(\sigma)\|_2d\sigma+
  \int_0^t \|\Delta \tilde u_\kappa(\sigma)\|_2\|\Psi_t(\sigma)\|_2d\sigma \\
  &\leq  
  \int_0^t \|\nabla\tilde u_\kappa(\sigma)\|_2\|\nabla\Psi_t(\sigma)\|_2d\sigma+
  \int_0^t \|\partial_t\tilde u_\kappa(\sigma)\|_2\|\Psi_t(\sigma)\|_2d\sigma \\
&+c\kappa\int_0^t  \|\Psi_t(\sigma)\|_2d\sigma \\
&\leq
\int_0^t \|\nabla\tilde u_\kappa(\sigma)\|_2\|\nabla\Psi_t(\sigma)\|_2d\sigma+
  \eps\int_0^t \|\partial_t\tilde u_\kappa(\sigma)\|_2^2d\sigma \\
 &+c_\eps\int_0^t\|\Psi_t(\sigma)\|_2^2d\sigma+
c\kappa\int_0^t  \|\Psi_t(\sigma)\|_2d\sigma,
\end{align*}
where $\Psi_t$ is the harmonic extension of $\psi_t$ to $\Omega$ 
(see formula \eqref{harm-ext-psi}). Analogously, we reach
\begin{align*}
\left|  \int_0^t 
{}_{H^{-\frac 12}\!\!}\bigg\langle
   \frac{\partial \tilde v_\kappa(\sigma)}{\partial\nu},\zeta_t(\sigma)
 \bigg\rangle_{\!\!H^{\frac 12}}\!\!d\sigma\right| &\leq
\int_0^t \|\nabla\tilde v_\kappa(\sigma)\|_2\|\nabla Z_t(\sigma)\|_2d\sigma+
  \eps\int_0^t \|\partial_t\tilde v_\kappa(\sigma)\|_2^2d\sigma \\
&+c_\eps\int_0^t\|Z_t(\sigma)\|_2^2d\sigma+
c\kappa\int_0^t  \|Z_t(\sigma)\|_2d\sigma,
\end{align*}
where, instead, $Z_t$ denotes the harmonic extension of $\zeta_t$ to $\Omega$, namely
\begin{equation*}
\begin{cases}
-\Delta Z_t(x;t)=0,&\text{in }\Omega,\\
\noalign{\vskip2pt}
\,\, Z_t(x;t)=\zeta_t(x,t), &\text{on }\partial\Omega.
\end{cases}
\end{equation*}
From the above estimates, the definition of 
$\Lambda_\kappa$, \eqref{eq:dinint}, Lemma \ref{stabilU}, 
and assumptions \eqref{eq:1} we obtain that
\begin{align*}
\|\nabla \tilde u_\kappa(t)\|_2^2 &+\|\nabla \tilde v_\kappa(t)\|_2^2\leq
C_1+C_2\int_0^t\|\nabla \tilde u_\kappa(\sigma)\|_{2}\big[\|\nabla
   U_t(\sigma)\|_{2}+\|\nabla \Psi_t(\sigma)\|_{2}\big]
d\sigma \\
&+C_3\int_0^t\|\nabla \tilde v_\kappa(\sigma)\|_{2}\big[\|\nabla
   V_t(\sigma)\|_{2}+\|\nabla Z_t(\sigma)\|_{2}\big]
d\sigma ,
\end{align*}
for some positive constant $C_1=C_1(\kappa)$ independent of $t$, 
\begin{equation}
\label{constantdip}
C_1(\kappa)=Q+\kappa\mu,\qquad Q=Q(u_0,v_0,\psi,\zeta),
\end{equation}
for  $C_2,C_3$ independent of $t$ and $\kappa$, where $\mu$ has been
defined in \eqref{cruquant}. Hence, by the Cauchy-Schwarz inequality
\begin{multline*}
\|\nabla \tilde u_\kappa(t)\|_2^2 +\|\nabla \tilde v_\kappa(t)\|_2^2\leq
C_1+C_4\int_0^t\sqrt{\|\nabla \tilde u_\kappa(\sigma)\|_2^2 
+\|\nabla \tilde v_\kappa(\sigma)\|_2^2 }\times\\
\times
\big[\|\nabla U_t(\sigma)\|_{2}+\|\nabla\Psi_t(\sigma)\|_{2}+\|\nabla
   V_t(\sigma)\|_{2}+\|\nabla Z_t(\sigma)\|_{2}\big]
\,d\sigma 
\end{multline*}
for all $t\geq 0$, for some positive constant  $C_4$ independent
of $t$ and $\kappa$. From Lemma \ref{gronwalltype} it follows that, for all $t\geq0$,
\begin{align*}
\|\nabla \tilde u_\kappa(t)\|_2^2 +\|\nabla \tilde v_\kappa(t)\|_2^2 &\leq 2C_1+C_4^2\big[ \|\nabla U_t\|_{L^1(0,\infty;L^2(\Omega))} 
 +\|\nabla\Psi_t\|_{L^1(0,\infty;L^2(\Omega))} \\
&+\|\nabla V_t\|_{L^1(0,\infty; L^2(\Omega))}+\|\nabla Z_t\|_{L^1(0,\infty;L^2(\Omega))}\big]^2,
\end{align*}
which, by Lemma \ref{stabilU1} and assumption \eqref{eq:1}, yields boundedness
of $(\tilde u_\kappa(t),\tilde v_\kappa(t))$ and, consequently 
of the sequence $(u_\kappa(t), v_\kappa(t)$ in $\Hh$, with the estimate
appearing in \eqref{acca1}.
\vskip2pt
In particular, from  the $\Hh$ boundedness of  $(\tilde u_\kappa(t),\tilde v_\kappa(t))$ we deduce that 
$\Lambda_\kappa$ is bounded from below uniformly with respect to $t$, with
a constant $\alpha_\kappa$ of the same form as the one appearing 
in inequality \eqref{acca1} (say, $\Lambda_\kappa(t)\geq -M-N\k\mu$, for some constants $M,N\geq 0$).
To prove this it suffices to repeat the estimates that we have obtained above (see the inequalities
following formula \eqref{constantupident}) on the
term which appear in the functional as time integrals, using the $H^1$
bound of $\tilde u_\kappa$ and $\tilde v_\kappa$, uniform in time. Notice that
the time integrals $\eps\int_0^t \|\partial_t\tilde u_\kappa(\sigma)\|_2^2d\sigma$ and
$\eps\int_0^t \|\partial_t\tilde v_\kappa(\sigma)\|_2^2d\sigma$ which appear
in the estimate of the boundary term are balanced by the corresponding
term in the definition of $\Lambda_\kappa$. 
More precisely, we obtain
$$
2\left|\int_0^t\int_{\Omega} \nabla \tilde u_\kappa(\sigma)\cdot\nabla U_t(\sigma)d\sigma\right|
\leq 2(R+\k\mu)\|\nabla U_t\|_{L^1(0,\infty;L^2(\Omega))}\leq A+B\k\mu,
$$
$$
\left|\int_\Omega \nabla U(t)\cdot\nabla \tilde u_\kappa(t)\right|\leq (R+\k\mu) \sup_{t\geq 1}\|\nabla U(t)\|_2
\leq C+D\k\mu,
$$
as well as
\begin{align*}
&-\int_0^t {}_{H^{-\frac 12}\!\!}\bigg\langle
   \frac{\partial \tilde u_\kappa(\sigma)}{\partial\nu},\psi_t(\sigma)
 \bigg\rangle_{\!\!H^{\frac 12}}\!d\sigma
 -\int_0^t {}_{H^{-\frac 12}\!\!}\bigg\langle
  \frac{\partial \tilde v_\kappa(\sigma)}{\partial\nu},\zeta_t(\sigma)
\bigg\rangle_{\!\!H^{\frac 12}}\!d\sigma  \\
&+\eps\int_0^t\|\partial_t \tilde u_\kappa(\sigma)\|_2^2\,d\sigma
+\eps\int_0^t\|\partial_t \tilde u_\kappa(\sigma)\|_2^2\,d\sigma \geq \\
& -(R+\k\mu) \|\nabla\Psi_t\|_{L^1(0,\infty, L^2(\Omega))}
 -c_\eps \|\Psi_t\|_{L^2(0,\infty, L^2(\Omega))}^2-
c\kappa \|\Psi_t\|_{L^1(0,\infty, L^2(\Omega))} \\
&-(R+\k\mu)\|\nabla Z_t\|_{L^1(0,\infty, L^2(\Omega))}
-c_\eps \|Z_t\|_{L^2(0,\infty, L^2(\Omega))}^2-
c\kappa \|Z_t\|_{L^1(0,\infty, L^2(\Omega))}\geq \\
& -E-F\k\mu,
\end{align*}
for some constants $A,B,C,D,E,F\geq 0$ independent of $\k$ and $t$.
Now, for all $\tau\in[0,\tau_0]$,
\begin{align*}
&\|\tilde u_\kappa^\tau(t)-\tilde u_\kappa(t)\|_2^2+\|\tilde v_\kappa^\tau(t)-\tilde v_\kappa(t)\|_2^2 \\
\noalign{\vskip4pt}
&=\int_\Omega|\tilde u_\kappa^\tau(t)-\tilde u_\kappa(t)|^2+
\int_\Omega|\tilde v_\kappa^\tau(t)-\tilde v_\kappa(t)|^2\\
&\leq \tau\int_t^{t+\tau}\|\partial_t \tilde u_\kappa(\sigma)\|^2_2\,d\sigma
+\tau\int_t^{t+\tau}\|\partial_t \tilde v_\kappa(\sigma)\|^2_2\,d\sigma \\
\noalign{\vskip2pt}
&= {\textstyle\frac{\tau}{1-\eps}}  \int_t^{t+\tau}\Big(-\frac{d}{d\sigma}\Lambda_\kappa(\sigma)\Big)\,d\sigma
\leq {\textstyle\frac{\tau_0}{1-\eps}}\big[{\Lambda_\kappa(t)-\Lambda_\kappa(t+\tau_0)}\big],
\end{align*}
where we exploited H\"older inequality, Fubini's Theorem and identity \eqref{ident} 
(in the spirit of \cite{cortazelgue}). Hence, we obtain
\begin{align}
\label{eq:16}
  &\|u_\kappa^\tau(t)- u_\kappa(t)\|_2^2+\| v_\kappa^\tau(t)- v_\kappa(t)\|_2^2 \\
  \noalign{\vskip5pt}
  &\leq2\big(\|\tilde u_\kappa^\tau(t)-
\tilde u_\kappa(t)\|_2^2+\|\tilde v_\kappa^\tau(t)-\tilde v_\kappa(t)\|_2^2+
\|U^\tau(t)-U(t)\|_2^2+\|V^\tau(t)-V(t)\|_2^2\big) \notag \\
  \noalign{\vskip5pt} &\leq {\textstyle\frac{2\tau_0}{1-\eps}}
  \big[{\Lambda_\kappa(t)-\Lambda_\kappa(t+\tau_0)}\big]+2\|U^\tau(t)-U(t)\|_2^2+2
  \|V^\tau(t)-V(t)\|_2^2.  \notag
\end{align}
Since $\Lambda_\kappa$ is nonincreasing and bounded from below at fixed $\k$, it follows that
$\Lambda_\kappa(t)$ admits a finite limit as $t\to\infty$. Therefore, letting $t\to\infty$ in \eqref{eq:16},
and taking into account Lemma \ref{stabilU}, we obtain \eqref{firstst}.
Now, assume by contradiction that, for some $\eps_0>0$,
$$
\|u_\kappa(t_h+\tau_h)-u_\kappa(t_h)\|_{H^1}\geq\eps_0>0,
$$
along a diverging sequence $(t_h)\subset\R^+$ and for $(\tau_h)\subset\R^+$ bounded.
In light of Lemma~\ref{compemb}, there exist $\hat u$ and $\check u \in H^1(\Omega)$ such that,
up to a subsequence that we still denote by $(t_h)$,
$u_\kappa(t_h+\tau_h)\to\hat u$ in $H^1(\Omega)$ as $h\to\infty$, and
$u_\kappa(t_h)\to\check u$ in $H^1(\Omega)$ as $h\to\infty$. In particular,
$\|\hat u-\check u\|_{H^1}\geq\eps_0>0,$
while \eqref{firstst} yields $\|\hat u-\check u\|_{L^2}=0$, thus
giving rise to a contradiction. One argues similarly for $v_\kappa$.
This concludes the proof of the theorem.
\end{proof}

\medskip
\noindent

Next we have an important consequence of the previous lemma, proving Theorem \ref{cve-cor}.

\begin{corollary}
\label{stabbzbis}
Assume  that $(u_0,v_0)\in\Hh_0$ and that the boundary data are stationary. 
Then the sequence $(u_\kappa(t),v_\kappa(t))$ is uniformly bounded in $H^1$ with respect to $t$ and $\kappa$.
Moreover the energy functional which appears in the proof of Theorem \ref{stabbz} is bounded below and above
by constants which are independent of $\kappa$.
\end{corollary}
\begin{proof}
If $(u_0,v_0)\in\Hh_0$, since $u_0v_0=0$ and $\psi_t=\zeta_t=0$ by \eqref{cruquant} we 
have that $\mu=0$. In turn, by  \eqref{acca1}, the sequence $(u_\kappa(t),v_\kappa(t))$
is uniformly bounded with respect to $t$ and $\kappa$. By inspecting the proof of Theorem \ref{stabbz}
it is easy to check that the auxiliary energy functional satisfies 
$$
-M-N\k\mu \leq \Lambda_\k(t)\leq O+P\k\mu, \quad t\geq 0,
$$
for some constants $M,N,O,P\geq 0$ independent of $\k$. Hence, being $\mu=0$
it follows that $\Lambda_\k$ has bounds uniform in time and in $k$.
\end{proof}

\subsection{Proof of Theorem  \ref{cve} concluded}
Let  $\kappa>0$ and let $(t_h)\subset\R^+$ be any diverging sequence.
Then, by virtue of Theorem~\ref{stabbz}, we have
\begin{equation}
\label{stabl}
\lim_{h\to\infty}\|u_\kappa(t_h+\tau_h)-u_\kappa(t_h)\|_{H^1}=0,\qquad
\lim_{h\to\infty}\|v_\kappa(t_h+\tau_h)-v_\kappa(t_h)\|_{H^1}=0,
\end{equation}
for every sequence $(\tau_h)\subset[0,1]$.
Let us fix $\eta,\xi\in C^\infty_c(\Omega)$. By integrating over $(t_h,t_h+1)\times\Omega$
the equations of  \eqref{system} multiplied by $\eta$ and $\xi$ respectively,  we reach
$$
\lim_h\left[\int_{t_h}^{t_h+1}\!\!\!\int_\Omega \partial_tu_\kappa\eta+
\int_{t_h}^{t_h+1}\!\!\!\int_\Omega \nabla u_\kappa\cdot\nabla\eta
-\int_{t_h}^{t_h+1}\!\!\!\int_\Omega f(u_\kappa)\eta+
\kappa\int_{t_h}^{t_h+1}\!\!\!\int_\Omega u_\kappa v_\kappa^2\eta\right]=0,
$$
$$
\lim_h\left[\int_{t_h}^{t_h+1}\!\!\!\int_\Omega \partial_tv_\kappa\xi+
\int_{t_h}^{t_h+1}\!\!\!\int_\Omega \nabla v_\kappa\cdot\nabla\xi
-\int_{t_h}^{t_h+1}\!\!\!\int_\Omega g(v_\kappa)\xi+
\kappa\int_{t_h}^{t_h+1}\!\!\!\int_\Omega v_\kappa u_\kappa^2\xi\right]=0.
$$
Regarding the first terms in the previous identities, we obtain
$$
\lim_{h}\left|\int_{t_h}^{t_h+1}\!\!\!\int_\Omega \partial_tu_\kappa\eta\right|\leq
\lim_{h}\int_\Omega |u_\kappa(t_h)-u_\kappa(t_h+1)||\eta|
\leq \lim_{h}c\|u_\kappa(t_h)-u_\kappa(t_h+1)\|_2=0,
$$
$$
\lim_{h}\left|\int_{t_h}^{t_h+1}\!\!\!\int_\Omega \partial_tv_\kappa\xi\right|\leq
\lim_{h}\int_\Omega |v_\kappa(t_h)-v_\kappa(t_h+1)||\xi|
\leq \lim_{h}c\|v_\kappa(t_h)-v_\kappa(t_h+1)\|_2=0.
$$
Moreover, there exist two sequences $(s_h),(r_h)\subset\R^+$ such that
$$
t_h\leq s_h\leq t_h+1,\quad t_h\leq r_h\leq t_h+1,\quad
s_h=t_h+\rho_h^1,\,\,\,r_h=t_h+\rho_h^2,
$$
with $(\rho_h^1),(\rho_h^2)\subset[0,1]$, and
$$
\int_{t_h}^{t_h+1}\!\!\int_\Omega \nabla u_\kappa\cdot\nabla\eta-
f(u_\kappa)\eta+\kappa u_\kappa v_\kappa^2\eta=\int_\Omega\nabla u_\kappa(s_h)\cdot\nabla\eta-
f(u_\kappa(s_h))\eta+\kappa u_\kappa(s_h) v_\kappa^2(s_h)\eta,
$$
$$
\int_{t_h}^{t_h+1}\!\!\int_\Omega\left[ \nabla v_\kappa\cdot\nabla\xi-
g(v_\kappa)\xi+\kappa v_\kappa u_\kappa^2\xi\right]=\int_\Omega\nabla v_\kappa(r_h)\cdot\nabla\xi-
g(v_\kappa(r_h))\xi+\kappa v_\kappa(r_h) u_\kappa^2(r_h)\xi.
$$
In turn, we get
$$
\lim_h\left[\int_\Omega \nabla u_\kappa(s_h)\cdot\nabla\eta
-\int_\Omega f(u_\kappa(s_h))\eta+
\kappa\int_\Omega u_\kappa(s_h) v_\kappa^2(s_h)\eta\right]=0,
$$
$$
\lim_h\left[\int_\Omega \nabla v_\kappa(r_h)\cdot\nabla\xi-
\int_\Omega g(v_\kappa(r_h))\xi+
\kappa\int_\Omega v_\kappa(r_h) u_\kappa^2(r_h)\xi\right]=0.
$$
On the other hand, in light of \eqref{stabl}, there holds
$$
\lim_h\|\nabla u_\kappa(s_h)-\nabla u_\kappa (t_h)\|_{2}
=\lim_h\|\nabla u_\kappa(t_h+\rho_h^1)-\nabla u_\kappa (t_h)\|_{2}=0,
$$
$$
\lim_h\|\nabla v_\kappa(r_h)-\nabla v_\kappa (t_h)\|_{2}
=\lim_h\|\nabla v_\kappa(t_h+\rho_h^2)-\nabla v_\kappa (t_h)\|_{2}=0.
$$
Hence $\int_\Omega (\nabla u_\kappa(s_h)-\nabla u_\kappa(t_h))\cdot\nabla\eta \to 0$ and, 
as $f,g$ are $C^1$ on $[0,1]$ and $0\leq u_\kappa,v_\kappa\leq 1$,
$$
\left|\int_\Omega  (f(u_\kappa(s_h))-f(u_\kappa(t_h)))\eta\right|\leq c\sup_{[0,1]} |f'|
\|u_\kappa(s_h)-u_\kappa(t_h)\|_2\to 0,
$$
as $h\to\infty$, and, finally,
$$
\left|\int_\Omega (u_\kappa(s_h) v_\kappa^2(s_h)-u_\kappa(t_h) v_\kappa^2(t_h))\eta\right|
\leq c\|u_\kappa(s_h)-u_\kappa(t_h)\|_2+c\|v_\kappa(s_h)-v_\kappa(t_h)\|_2\to 0,
$$
as $h\to\infty$, the positive constant $c$ varying from line to line. 
Of course, the same conclusions hold for the limit involving 
the sequence $v_\kappa(r_h)$. In conclusion, we reach
$$
\lim_h\left[\int_\Omega \nabla u_\kappa(t_h)\cdot\nabla\eta
-\int_\Omega f(u_\kappa(t_h))\eta+
\kappa\int_\Omega u_\kappa(t_h) v_\kappa^2(t_h)\eta\right]=0,
$$
$$
\lim_h\left[\int_\Omega \nabla v_\kappa(t_h)\cdot\nabla\xi-
\int_\Omega g(v_\kappa(t_h))\xi+
\kappa\int_\Omega v_\kappa(t_h) u_\kappa^2(t_h)\xi\right]=0.
$$
Again in view of Theorem \ref{stabbz}, we can assume that, up to a
subsequence, which we shall denote again by $t_h$, we have that
$u_\kappa(t_h)\rightharpoonup \hat u_\kappa$ and
$v_\kappa(t_h)\rightharpoonup \hat v_\kappa$ weakly in $H^1(\Omega)$.
Up to a subsequence, in light of Lemma \ref{compemb}, this convergence
is actually strong. Notice also that
\begin{align*}
& \hat u_\kappa|_{\partial\Omega}=\lim_h u_\kappa(t_h)|_{\partial\Omega}=\lim_h \psi(t_h)|_{\partial\Omega}=\psi_\infty, \\
& \hat v_\kappa|_{\partial\Omega}=\lim_h v_\kappa(t_h)|_{\partial\Omega}=\lim_h \zeta(t_h)|_{\partial\Omega}=\zeta_\infty,
\end{align*}
where we exploited the compact embedding $H^1(\Omega)\hookrightarrow H^{1/2}(\partial\Omega)$.
Moreover, by Lemma \ref{l:mp} and the Dominated Convergence Theorem,
as $h\to \infty$, we get 
$$
\int_\Omega \nabla \hat u_\kappa\cdot\nabla\eta
-\int_\Omega f(\hat u_\kappa)\eta+
\kappa\int_\Omega \hat u_\kappa \hat v_\kappa^2\eta=0,\qquad\forall\eta\in H^1_0(\Omega),
$$
$$
\int_\Omega \nabla \hat v_\kappa\cdot\nabla\xi-
\int_\Omega g(\hat v_\kappa)\xi+
\kappa\int_\Omega \hat v_\kappa \hat u_\kappa^2\xi=0,\qquad\forall\xi\in H^1_0(\Omega).
$$
Hence $(\hat u_\kappa,\hat v_\kappa)\in \Hh$ is a solution to \eqref{limsystem}.  The convergence
occurs of course in $L^p(\Omega)$ for any $p\in[2,2^*)$. For $p\geq 2^*$,  
taking $\eps>0$ and using the bounds 
$0\leq u_\kappa(t_h) \leq 1$ and $0\leq \hat u_\k \leq 1$, we have
$$
\int_\Omega |u_\kappa(t_h)-\hat u_\k|^p\leq 
2^{p+\eps-2^*}\|u_\kappa(t_h)-\hat u_\k\|_{2^*-\eps}^{2^*-\eps},
$$
concluding the proof.
\qed

\bigskip

\section{Proof of Theorem~\ref{main1}}

Before concluding the proof of Theorem ~\ref{main1}, we provide the  convergence of the
sequences $(\hat u_\kappa,\hat v_\kappa)$ in any $L^p$ space with $p\geq 2$  towards a segregated state.
Notice that the solutions $(\hat u_\kappa,\hat v_\kappa)$ to \eqref{limsystem}
pop up as $H^1$ limits of the solutions to \eqref{system-BDtime ind}, and the 
boundedness  of $(\hat u_\kappa,\hat v_\kappa)$ in $H^1$ is inherited by the 
boundedness of $(u_\kappa(t_h),v_\kappa(t_h))$ in $H^1$ uniform in $t$ and $\kappa$ (in the case $(u_0,v_0)\in\Hh_0$).
Without this information it would not have been clear how to show the boundedness of $(\hat u_\kappa,\hat v_\kappa)$ 
working directly on the elliptic system (instead, for system \eqref{orig2po}, this is an
easy task, cf.\ \cite[Lemma 2.1]{ctv-adv}). 

\begin{lemma}
\label{cvell}
Assume that $(u_0,v_0)\in\Hh_0$. Let $(\hat u_\kappa,\hat v_\kappa)\in\Hh$ be the solution
to system \eqref{limsystem} as obtained in Theorem \ref{cve} for $\k>0$. Then
there exists $(u_\infty,v_\infty)\in\Hh_0$ with 
$$
u_\infty,v_\infty\geq 0,\quad  -\Delta u_\infty\leq f(u_\infty),\quad -\Delta v_\infty\leq g(v_\infty)
$$
and $u_\infty|_{\partial\Omega}=\psi$, $v_\infty|_{\partial\Omega}=\zeta$ 
such that, up to a subsequence, as $\kappa\to\infty$,
$$
(\hat u_\kappa,\hat v_\kappa)\to (u_\infty,v_\infty)
\quad \text{in the $L^p\times L^p$ norm
for any $p\in[2,\infty)$.}
$$
\end{lemma}

\begin{proof}
By virtue of Corollary \ref{stabbzbis} the sequence $(u_\kappa(t_h),v_\kappa(t_h))$ is bounded in
$\Hh$, uniformly with respect to $\kappa$. Hence, since $(\hat u_\kappa, \hat v_\kappa)$ is the 
$H^1$-limit of $(u_\kappa(t_h),v_\kappa(t_h))$ as $h\to\infty$, we deduce that $(\hat u_\kappa,\hat v_\kappa)$ is
bounded in $\Hh$ and $0\leq \hat u_\k(x)\leq 1$,  $0\leq \hat v_\k(x)\leq 1$, for a.e.\ $x\in\Omega$. 
Taking into account that some terms in the functional $\Lambda_\k$ introduced within the proof of 
Theorem \ref{stabbz} vanish under the current assumptions 
(stationary boundary conditions) and that the terms $\eps\int_0^t\|\partial_t \tilde u_\kappa(\sigma)\|_2^2$
and $\eps\int_0^t\|\partial_t \tilde v_\kappa(\sigma)\|_2^2$ were artificially attached 
to make things work (notice that the original $\Lambda_\k$ 
is decreasing also in the case $\eps=0$, see formula \eqref{ident}), we now just consider
the natural energy functional (for the sake of simplicity we do not change the name)
\begin{align*}
  \Lambda_\kappa(t) &=\frac 12\|\nabla \tilde u_\kappa(t)\|_2^2 +\frac
  12\|\nabla \tilde v_\kappa(t)\|_2^2-
  \int_\Omega F(\tilde u_\kappa(t)+U(t))\\
  &\,\,\, -\int_\Omega G(\tilde v_\kappa(t)+V(t))+
  \frac \kappa2\int_\Omega (\tilde u_\kappa(t)+U(t))^2(\tilde v_\kappa(t)+V(t))^2.
\end{align*}
Then, we have
$$
\kappa\int_\Omega u_\kappa^2(t_h)v_\kappa^2(t_h)=2\Lambda_\kappa(t_h)-
\|\nabla \tilde u_\kappa(t_h)\|_2^2-\|\nabla \tilde v_\kappa(t_h)\|_2^2+
2\int_\Omega F(u_\kappa(t_h))+G(v_\kappa(t_h)).
$$
Since by Corollary \ref{stabbzbis} the right hand side is uniformly 
bounded with respect to $\kappa$, we have
\begin{equation}
\label{eq:20}
\kappa\int_\Omega \hat u_\kappa^2 \hat v_\kappa^2\leq C,
\end{equation}
for some positive constant $C$ independent of $\kappa$.
Let $u_\infty\in H^1(\Omega)$ and $v_\infty\in H^1(\Omega)$ be
the weak limits, as $\kappa\to\infty$, of $\hat u_\kappa$ and $\hat v_\kappa$ in $H^1(\Omega)$ respectively.
Or course, by the compact embedding $H^1(\Omega)\hookrightarrow L^{2^*}(\Omega)$, up to a further 
subsequence, $\hat u_\kappa\to u_\infty$ and $\hat v_\kappa\to v_\infty$ in $L^p(\Omega)$
for any $p\in[2,2*)$ and $0\leq u_\infty(x)\leq 1$,  $0\leq v_\infty(x)\leq 1$, for a.e.\ $x\in\Omega$. 
In the case $p\geq 2^*$, let $\eps>0$, so that
$$
\int_\Omega |\hat u_\kappa-u_\infty|^p\leq 2^{p+\eps-2^*}\|\hat u_\kappa-u_\infty\|_{2^*-\eps}^{2^*-\eps},
$$
yielding again the convergence. Due to inequality \eqref{eq:20}, we get
$$
\lim_{\kappa\to\infty}\int_\Omega \hat u_\kappa^2\hat v_\kappa^2=\int_\Omega u_\infty^2v_\infty^2=0,
$$
which yields $u_\infty v_\infty=0$ a.e.\ in $\Omega$, namely $(u_\infty,v_\infty)\in\Hh_0$. 
Moreover, for each $\kappa>0$, 
\begin{equation*}
-\Delta \hat u_\k\leq f(\hat u_\k),\qquad
-\Delta \hat v_\k\leq g(\hat v_\k),
\end{equation*}
which pass to the weak the limit, yielding $-\Delta u_\infty\leq f(u_\infty)$ 
and $-\Delta v_\infty\leq g(v_\infty)$.  By the compact embedding 
$H^1(\Omega)\hookrightarrow H^{1/2}(\partial\Omega)$, also
the boundary conditions are conserved.
\end{proof}

\subsection{Proof of Theorem~\ref{main1} concluded}
We can now conclude the proof of Theorem~\ref{main1}.
Let $(u_0,v_0)\in\Hh_0$, $p\in[2,\infty)$ and let $(t_h)\subset\R^+$ be any diverging
sequence. In light of Theorem \ref{cve}, for every $\kappa\geq 1$, 
there exist a solution $(\hat u_\kappa,\hat v_\kappa)$ of \eqref{limsystem}
and a subsequence $(t^\kappa_h)\subset\R^+$ such that,
$$
\|(u_{\kappa}(t_h^\kappa),v_\kappa(t_h^\kappa)) -(\hat u_\kappa,
\hat v_\kappa)\|_{\Hh} \to 0,\quad\text{as $h\to\infty$}.
$$
Moreover, by Lemma \ref{cvell}, there exists $(u_\infty,v_\infty)\in\Hh_0$ with the required
properties, such that, up to a subsequence, 
$$
\|(\hat u_\kappa,\hat v_\kappa)-(u_\infty,v_\infty)\|_{L^p\times L^p}\to 0,\quad\text{as $\kappa\to\infty$}.
$$
Now, let $m\geq 1$ and let $\kappa_m\geq 1$ be such that
$$
\|(\hat u_{\kappa_m},\hat v_{\kappa_m})-(u_\infty,v_\infty)\|_{L^p\times L^p}<\frac{1}{2m}.
$$
Then, there exists $t_{h_m}^{\kappa_m}\geq 1$ such that
$$
\|(u_{\kappa_m}(t_{h_m}^{\kappa_m}),v_{\kappa_m}(t_{h_m}^{\kappa_m})) -
(\hat u_{\kappa_m},\hat v_{\kappa_m})\|_{L^p\times L^p}<\frac{1}{2m}.
$$
In turn, setting $t_m=t_{h_m}^{\kappa_m},$
and combining the previous inequalities, we get
$$
\|(u_{\kappa_m}(t_m),v_{\kappa_m}(t_m)) -(u_\infty,v_\infty)\|_{L^p\times L^p}<\frac{1}{m},
$$
which concludes the proof of the first assertion.  In the one dimensional
case, by means of Morrey Theorem, for every $x,y\in\Omega$, we have
$$
|u_{\kappa_m}(t_m)(x)-u_{\kappa_m}(t_m)(y)|
\leq 4\|\nabla u_{\kappa_m}(t_m)\|_2\sqrt{|x-y|}\leq C\sqrt{|x-y|},
$$
together with $|u_{\kappa_m}(t_m)(x)|\leq 1$, yielding the 
convergence to $(u_\infty,v_\infty)$ in the $L^\infty\times L^\infty$ 
norm via Ascoli's Theorem.
\qed

\bigskip
\subsection*{Acknowledgment}
The author is indebted with the anonymous Referee for a very careful reading of the 
manuscript and for many valuable suggestions and comments which helped to improve
the paper.

The author wishes to thank Prof.\ Yoshio Yamada for some useful references concerning 
the existence of global solutions to reaction diffusion systems and Prof.\ Alain Haraux 
for some comments about the long term behaviour in presence of time-dependent boundary data.

\bigskip

\end{document}